\documentclass[12pt]{article}

\usepackage{color}
\usepackage[T2A]{fontenc}
\usepackage[english, russian]{babel}
\usepackage{amsmath}
\usepackage{amssymb}
\usepackage{amsthm}
\usepackage{dsfont}
\usepackage{textcomp}

\begin{document}

\title {Множественная сложность построения правильного многоугольника
\footnote{ 
Благодарю Д. Мусатова, А. Скопенкова и А. Савватеева за ценные замечания и предложения при написании данной работы
}
}

\author{Е.С. Коган}
\date{}

\maketitle

\begin{abstract}
Given a subset of $\mathds C$ containing $x,y$,  one can add $x + y,\,x - y,\,xy$ or (when $y\ne0$) $x/y$
or any $z$ such that $z^2=x$.
Let $p$ be a prime Fermat number.
We prove that it is possible to obtain from $\{1\}$ a set containing all the $p$-th roots of 1
by $12 p^2$ above operations.
This result is different from the standard estimation of complexity of an algorithm computing the $p$-th roots of 1.
\end{abstract}


К подмножеству $A\subset \mathds C$ содержащему числа $x,y$ можно добавить любое из чисел $x + y,\,x - y,\,xy$ или (если $y\ne0$) $x/y$
или любое $z$ такое, что $z^2=x$.


\smallskip

\textbf{Основная теорема.}
Пусть $p$ --- простое число Ферма, т.е. простое число вида $2^m + 1$, где $m$ --- степень 2,
$\varepsilon := \cos\frac{2\pi}{p} + i \sin\frac{2\pi}{p}$.
Тогда из $\{1\}$ можно получить некоторое множество, содержащее числа $1, \varepsilon, \varepsilon^2, \ldots, \varepsilon^{p - 1}$,
за $12 p^2$ добавлений, определенных выше.

\bigskip



Эта работа может быть интересна читателю, так как на примере решения алгоритмической задачи иллюстрирует важный математический метод.

Назовем \emph{множественной сложностью} сложность, рассмотренную в основной теореме.
Такое понятие сложности отличается от сложности как времени работы алгоритма, находящего корни степени $p$ из 1. Однако последняя сложность также пропорциональна $p^2$. Об алгоритмах вычисления корней $p$-й степени из 1
см. \cite{safin} и, возможно, \cite{algorithm}. О строгих определениях различных понятий сложности см. \cite{complexity}. Видимо, понятие множественной сложности совпадает с одним из сформулированных в этой книге, возможно, со сложностью как глубиной формулы. В любом случае, основная теорема не претендует на новизну. Даже если рассмотренное понятие сложности новое, доказательство не содержит новых идей.

\bigskip

\textbf{Замечание.}
Из основной теоремы можно вывести следующее утверждение:

Пусть $p$ --- простое число Ферма.
Тогда существует действительное число $C$, не зависящее от $p$, такое, что
за $C \cdot p^2$ операций проведения окружности с центром в одной точке и
проходящей через другую и проведения прямой через две точки из единичного
отрезка можно
получить правильный $p$-угольник.

Формализация аналогична основной теореме.

Известна история об аспиранте, который
разработал построение правильного многоугольника с $65537$ сторонами
за $20$ лет (см.~\cite{65537}).


\smallskip

\textbf{Замечание.}
Из основной теоремы также можно вывести ее вещественный аналог,
который состоит в следующем:

Существует такое число $C$, что для любого простого числа Ферма $p$
число $\cos{\frac{2\pi}{p}}$
можно получить из $\{1\}$ за $C \cdot p^2$ операций, аналогичным определенным выше, но при которых корень можно брать только
из положительных чисел.

\smallskip

Доказательство основной теоремы аналогично \cite[п.\,5.3.4]{source}.
Оценка же, получающаяся из доказательства построимости, приведенном в
\cite[конец п.\,5.3.3]{source}, пропорциональна $p^3$.

Следующее изложение идеи доказательства заимствованно из \cite[п. 5.3.4]{source}.

\smallskip

\emph{Идея доказательства основной теоремы для $p=5$.}
Сразу выразить число $\varepsilon$ через радикалы
трудно, поэтому сначала выразим некоторые <<многочлены от
$\varepsilon$>>. Мы знаем, что
$\varepsilon+\varepsilon^2+\varepsilon^3+\varepsilon^4=-1$. Поэтому
 $$
 (\varepsilon+\varepsilon^4)(\varepsilon^2+\varepsilon^3)=\varepsilon+\varepsilon^2+\varepsilon^3+\varepsilon^4=-1.
 $$
Обозначим
 $$
 T_0:=\varepsilon+\varepsilon^4\quad \text{и}\quad
T_1:=\varepsilon^2+\varepsilon^3.
 $$
Тогда по теореме Виета числа $T_0$ и~$T_1$ являются корнями уравнения
$t^2+t-1=0$. Поэтому можно выразить $T_0$ (и~$T_1$). Поскольку
$\varepsilon\cdot \varepsilon^4=1$, по теореме Виета числа
$\varepsilon$ и~$\varepsilon^4$ являются корнями уравнения
$t^2-T_0t+1=0$. Поэтому можно выразить $\varepsilon$
(и~$\varepsilon^4$).

\smallskip
\emph{Идея доказательства основной теоремы в общем случае.}
Сначала хорошо бы разбить сумму
 $$
 \varepsilon+\varepsilon^2+\ldots+\varepsilon^{p-1}=-1
 $$
на 2 слагаемых $T_0,T_1$, которые можно получить добавлениями, описанными перед теоремой 1 (иными словами, \emph{сгруппировать} хитрым образом корни
уравнения $1+x+x^2+\ldots+x^{p-1}=0$). Затем нужно разбить каждую
сумму $T_k$ на $2$ слагаемых $T_{k,0},T_{k,1}$,
которые можно получить такими добавлениями. И~так далее, пока не получим
$T_{\text{\scriptsize$\underbrace{1,\ldots,1}_s$}}=\varepsilon$.

\smallskip
Однако придумать нужные группировки чисел
$1,\varepsilon,\varepsilon^2,\ldots,\varepsilon^{p-1}$ "---
нетривиальная задача.

Теорема о~первообразном корне, приведенная, например, в \cite[п. 5.3.3]{source},
позволяет закодировать ненулевые вычеты по модулю $p$
вычетами по модулю $p-1$.
А~именно, выбрав первообразный корень $g$, мы вычету $k$ по модулю
$p-1$ сопоставляем (ненулевой) остаток от деления $g^k$ на~$p$.
Это кодирование фактически было использовано в~группировках,
построенных выше для $p=5$.

\bigskip

Введем определения и обозначения, необходимые для доказательства.

Множество $A \subset \mathds C$
\emph{построимо за $n$ операций из $B \subset \mathds C$},
если какое-нибудь множество $A' \supset A$ можно получить из $B$ за
$n$ операций, описанных выше.

Обозначим через

\begin{itemize}
\item
$g$ --- первообразный корень по модулю $p$\,;
\item
$T_{k, r} := \sum\limits_{a = 0}^{2^{m - k} - 1}{\varepsilon^{g^{2^k \cdot a+ r}}}$
для каждого $k \in \{0, 1, \ldots, m\}$, $r \in \mathds{Z}_{2^k}$.

В частности, $T_{0, 0} = -1$, а $T_{m, r} = \varepsilon^r$;
\item
$N_{k, t}$ --- число пар $(c, d)$ вычетов
по модулю $2^{m - k - 1}$, являющихся решениями сравнения
\begin{equation}
g^{2^{k + 1} \cdot c + t} + g^{2^{k + 1} \cdot d + 2^k + t}
\equiv 1 \pmod{p},
\label{raw_alpha}
\end{equation}
где $k \in \{0, 1, \ldots, m - 1\}$, $t \in \mathds{Z}_{2^m}$.
$N_{k, t}$ зависит от $m$, но поскольку $m$ зафиксировано, оно пропускается.

\end{itemize}

\textbf{Пример для $p=5$.}
Для любых вычетов $t_1, t_2$ по модулю 4, сравнимых по модулю 2, 
$N_{1, t_1} = N_{1, t_2}$.

Действительно, пусть $t_1 \equiv t_2 \pmod{4}.$
Тогда можно сопоставить каждому решению $(c, d)$ сравнения (\ref{raw_alpha}) при $k = 1,\ t = t_1$ решение $(c', d')$ сравнения (\ref{raw_alpha}) при $k = 1,\ t = t_2$ следующим образом.
Обозначим через $l~:=~(t_1 - t_2)/4,\ c'~:=~c~+~l,\ d'~:=~d~+~l$. Тогда:
\begin{equation*}
g^{4 c + t_1} + g^{4 d + 2 + t_1}
\equiv 1 \pmod{5}
\end{equation*}
\begin{equation*}
g^{4 (c + l) + t_2}
+ g^{4 (d + l) + 2 + t_2}
\equiv 1 \pmod{5}
\end{equation*}
\begin{equation*}
g^{4 c' + t_2}
+ g^{4 d' + 2 + t_2}
\equiv 1 \pmod{5}
\end{equation*}

Очевидно, что это сопоставление --- биекция. Cледовательно, $N_{1, t_1} = \\ = N_{1, t_2}$.

Осталось показать, 
что $N_{1, t} = N_{1, 2 + t}$.
Обозначим через $c' := d, d' := \\ := c - 1$.
Тогда:
\begin{equation*}
g^{4 c + t} + g^{4 d + 2 + t}
\equiv 1 \pmod{5}
\end{equation*}
\begin{equation*}
g^{4 d + (2 + t)}
+ g^{4 (c - 1) + 2 + (2 + t)}
\equiv 1 \pmod{5}
\end{equation*}
\begin{equation*}
g^{4 c' + (2 + t)}
+ g^{4 d' + 2 + (2 + t)}
\equiv 1 \pmod{5}
\end{equation*}

То есть мы, как и раньше, построили биекцию между решениями соответствующих сравнений.

\smallskip

\textbf{Лемма 1.}
Для любых вычетов $t_1, t_2 \in \mathds{Z}_{2^m}$, сравнимых по модулю $2^k$, 
$k \in \{0, 1, \ldots, m - 1\}$ верно
$N_{k, t_1} = N_{k, t_2}$.

\textbf{Доказательство.}
Пусть $t_1 \equiv t_2 \pmod{2^{k + 1}}.$
Обозначим через $l := \\ := \cfrac{t_1 - t_2}{2^{k + 1}}, c' := c + l, d' := d + l$.
Тогда:
\begin{equation*}
g^{2^{k + 1} \cdot c + t_1} + g^{2^{k + 1} \cdot d + 2^k + t_1}
\equiv 1 \pmod{p}
\end{equation*}
\begin{equation*}
g^{2^{k + 1} \cdot (c + l) + t_2}
+ g^{2^{k + 1} \cdot (d + l) + 2^k + t_2}
\equiv 1 \pmod{p}
\end{equation*}
\begin{equation*}
g^{2^{k + 1} \cdot c' + t_2}
+ g^{2^{k + 1} \cdot d' + 2^k + t_2}
\equiv 1 \pmod{p}
\end{equation*}

Cледовательно, $N_{k, t_1} = N_{k, t_2} $.

Для завершения доказательства леммы 1 достаточно показать,
что $N_{k, t} = N_{k, 2^k + t}$.
Обозначим через $c' := d, d' := c - 1$.
Тогда:
\begin{equation*}
g^{2^{k + 1} \cdot c + t} + g^{2^{k + 1} \cdot d + 2^k + t}
\equiv 1 \pmod{p}
\end{equation*}
\begin{equation*}
g^{2^{k + 1} \cdot d + (2^k + t)}
+ g^{2^{k + 1} \cdot (c - 1) + 2^k + (2^k + t)}
\equiv 1 \pmod{p}
\end{equation*}
\begin{equation*}
g^{2^{k + 1} \cdot c' + (2^k + t)}
+ g^{2^{k + 1} \cdot d' + 2^k + (2^k + t)}
\equiv 1 \pmod{p}
\end{equation*}
\qed

\bigskip

\textbf{Лемма 2.}
Для любого целого числа~$k$ от $1$ до $m - 1$ множество
$$A = \{0, 1,\,\dots\,, p\} \cup \{T_{k + 1, r}\ |\ r\,\in\,\mathds{Z}_{2^{k + 1}}\}$$
построимо из множества
$$B = \{0, 1,\,\dots\,, p\} \cup \{T_{k, r}\ |\ r\,\in\,\mathds{Z}_{2^k}\}$$
за $11 \cdot 4^k$ операций.

\smallskip

\textbf{Доказательство.}
Во-первых, для любых 
$k \in \{0, 1, \ldots, m\}$ и $t \in \mathds{Z}_{2^m}$
выполняется $N_{k, t} \leqslant p,$
т.\,к. в сравнении~{(\ref{raw_alpha})}
одному вычету $c$ может подходить не больше одного вычета $d.$
Следовательно, все $N_{k, t}$ содержатся в $B$.

Во-вторых,
\begin{equation*}
T_{k + 1, r} T_{k + 1, 2^k + r}
= \sum_{s = 0}^{2^m - 1} {N_{k, r - s} \varepsilon^{g^s}}
= \sum_{s = 0}^{2^k - 1} {N_{k, r - s} T_{k, s}}.
\end{equation*}
Здесь мы воспользовались леммой $1$.
При этом множество\\ $\{N_{k, r - s} T_{k, s}\ |\ r, s\in\mathds{Z}_{2^k}\}$ построимо из $B$ за $2^k \cdot 2^k$ операций, а значит, множество
$P := \{T_{k + 1, r} T_{k + 1, 2^k + r}\ |\ r\in\mathds{Z}_{2^k}\}$
построимо из $B$ за $2^k \cdot 2^k + + (2^k - 1) \cdot 2^k < 2 \cdot 4^k$ операций.

Далее, множество $P\cup B$ содержит
$$T_{k + 1, r}+\ T_{k + 1, 2^k + r} = T_{k, r} \in B\quad\text{и}\quad T_{k + 1, r} T_{k + 1, 2^k + r} \in P,$$
и для любых комплексных чисел $x_1, x_2$
множество $\{x_1, x_2\}$ построимо за $9$ операций из
множества $\{x_1 + x_2, x_1 x_2\}$
(по формуле $x_{1, 2}  = \frac{-b \pm \sqrt{b^2 - 4c}}{2}$).
Следовательно,
$\bigl\{T_{k + 1, r}\ |\ r\in \mathds{Z}_{2^{k+1}}\bigr\}$
построимо из $B$ за $2 \cdot 4^k + 9 \cdot 2^k \leqslant 11 \cdot 4^k$
операций. Кроме того, $\{0, 1,\,\dots\,,p\} \subset B$, поэтому
$A$ построимо из $B$ за $11 \cdot 4^k$ операций.
\qed

\smallskip

\textbf{Доказательство основной теоремы.}
Из леммы 2 следует, что множество
$$\bigl\{T_{m, r}\ |\ r \in \{0, 1, \ldots, 2^m - 1\}\bigr\} = \bigl\{\varepsilon^r\ |\ r \in \{0, 1, \ldots, 2^m - 1\}\bigr\}$$
построимо из $\{1\}$ за
\begin{equation*}
p + \sum_{k = 0}^{m - 1} {11 \cdot 4^k}
= p + 11 \frac{4^m - 1}{4 - 1}
< p + 11 \cdot 4^m < 12 p^2
\end{equation*}
операций.
\qed

\smallskip






\vspace*{\fill}

\noindent\rule{4cm}{0.5pt}

\bigskip

\noindentЕ.С. Коган, ГБОУ "Интеллектуал"

\noindent\texttt{koganeser@gmail.com}

\end{document}